\documentclass[oneside,11pt]{article} 

\usepackage{epsfig,graphicx}
\usepackage{latexsym,amssymb}
\usepackage{setspace,cite} 

\usepackage{amsmath, amssymb, amsthm}
\usepackage{graphicx,color}

\usepackage{fancyhdr}

\usepackage{a4}
\usepackage{bbm}

\newcommand{\R}{{^*\B{R}^{\textrm{nst}}}}
\newcommand{\B}[1]{\mathbb{#1}}
\usepackage{amsmath}
\usepackage{amssymb}
\usepackage{euscript}
\usepackage{amsthm}
\usepackage[all]{xy}
\newtheorem{prop}{Proposition}[section]

\newtheorem{cor}[prop]{Corollary}
\newtheorem{lemma}[prop]{Lemma}
\newtheorem{thm}[prop]{Theorem}

\theoremstyle{definition}
\newtheorem{defn}[prop]{Definition}

\theoremstyle{definition}
\newtheorem{rmk}[prop]{Remark}
\theoremstyle{definition}

\theoremstyle{definition}

\title{A Nonstandard Approach to Real Multiplication}
\author{Lawrence Taylor}

\begin{document}
\maketitle

\begin{abstract}
This paper looks at the use of Model Theoretic ideas to study certain non-Hausdorff spaces with a view to their application in Manin's proposed theory of Real Multiplication, which seeks to establish a framework in which to prove Hilbert's twelfth problem for real quadratic fields.  We study morphisms between these non-Hausdorff spaces using Nonstandard Analysis, and as a consequece in certain special cases we are able to describe an action of a certain Galois group on isomorphism classes of these spaces.
\end{abstract}

\section{Introduction}

Class Field Theory attempts to classify the abelian extensions of a field $K$, in terms of data intrinsic to $K$, namely the idele class group.  When $K$ is a number field, the celebrated ``Existence Theorem'' asserts a bijection between the finite abelian extensions over $K$, and open subgroups of finite index in the idele class group.  Fields arising in this way are known as \emph{class fields} for $K$.  Despite ensuring the existence of such fields for $K$, the proof of this result is in general nonconstructive, not providing a set of generators for the abelian extension.  In the early twentieth century Hilbert listed the need for such an \emph{explicit class field theory} as the twelfth in a list of twenty three problems he deemed to be of importance to mathematics.\newline

A century later Hilbert's twelfth problem remains unanswered, except in a few special circumstances.  In 1896 Hilbert himself gave the first complete answer to the case when $K$ is the field $\B{Q}$ of rational numbers following the work of Kronecker and Weber.  By the end of the nineteenth century a solution was known for the case when $K$ is an imaginary quadratic field, fulfilling Kronecker's \emph{Jugendtraum}, or ``dream of youth''.  This was achieved by generating abelian extensions of $K$ by adjoining special values of certain functions on elliptic curves with Complex Multiplication.  Much later in the 1970's, this result inspired the construction of the Lubin-Tate formal group, which was used to establish a solution when $K$ is a local field \cite{Iwasawa}.\newline

The simplest class of number fields for which the problem remains unsolved is the case when $K$ is a real quadratic field.\newline

Some progress has been made in obtaining solutions for isolated classes of real quadratic fields by Shimura \cite{ShimuraI,ShimuraII} and Shintani \cite{ShintaniIII}.  Shimura uses a similar philosophy to that of the theory of Complex Multiplication, generating abelian extensions of certain real quadratic fields (of class number one) by the torsion points of certain abelian varieties.  Shintani's work is motivated by the ability to express the L-function of a real quadratic field in terms of certain special functions studied by Barnes in \cite{BarnesII}.  However, neither of these provides a systematic way of giving a solution for a general real quadratic field.  \\

The quest for such a solution is the subject of Manin's paper ``Real Multiplication and Noncommutative Geometry'' \cite{Manin} where he poses his \emph{Alterstraum} - a theory of Real Multiplication.  The theory laid out in Manin's paper is connected to the development of Noncommutative Geometry studied by A.Connes in the early 1980's.  His book \cite{Connes}, considered one of the milestones in mathematics, studied the analytical theory of non-Hausdorff spaces using C*-algebra and operator theory.  Such ``Noncommutative spaces'' have a noncommutative C*-algebra associated to them which is an analogue to the algebra of $\B{C}$-valued functions on a Hausdorff space. \\ 

Whereas elliptic curves with Complex Multiplication play a key role in the solution of Hilbert's problem for quadratic imaginary fields, Manin's approach suggests the use of Noncommutative Tori in the analogous theory for real quadratic fields.  These are noncommutative C*-algebras parameterised by a real number $\theta$, generated by unitaries $U$ and $V$ satisfying the commutativity relation
$$VU=e^{2 \pi i \theta}UV.$$
We denote such an algebra by $A_\theta$.  Through the work of Reiffel and others a duality is established through such objects and isomorphism classes of pseudolattices, analogous to the relationship between elliptic curves over $\B{C}$ and lattices established by the Uniformisation Theorem.\\

The category of pseudolattices is defined to be
\begin{defn}[Category of Pseudolattices]\label{PL}
Let $\mathcal{PL}$ denote the category such that:
\begin{itemize}
\item The objects of $\mathcal{PL}$ are dense additive subgroups $L$ of $\B{R}$ of rank two.
\item A morphism between two pseudolattices $L_1$ and $L_2$ is a nonzero positive real number $\beta$
such that $\beta L_1 \subseteq L_2$.
\end{itemize}
\end{defn} 
The relationship between Noncommutative Tori and pseudolattices in Manin's paper is based on the association
$$A_\theta \leftrightarrow L_\theta:=\B{Z}+\theta\B{Z}.$$ 
Motivated by the duality between lattices and elliptic curves we make the following definition:

\begin{defn}[Quantum Torus]
Let $L$ be a pseudolattice.  The Quantum Torus associated to $L$ is the topological space $Z_L$ defined by $\B{R}/L$.
\end{defn} 

The non-Hausdorff nature of Quantum Tori imply that the space of continuous $\B{C}$-valued functions on $Z_L$ is trivial.  In \cite{FesenkoII} and \cite{FesenkoI}, Fesenko suggests the use of Model Theory, and more specifically Nonstandard Analysis as a tool for studying such topological spaces.  \newline

\subsection{Statement of results and overview}
In this paper we explore the use of nonstandard analysis to further our understanding of Quantum Tori from both a topological and arithmetic viewpoint.  Our approach is to consider a class of nonstandard topological spaces $T_L$ parameterised by pseudolattices, together with a surjective morphism whose image is a Quantum Torus.  These nonstandard topological spaces are called \emph{Hyper Quantum Tori}.\\

In \S\ref{monadicspace} we define the category $\mathcal{LIQ}$ of \emph{locally internal quotient spaces}.  Hyper Quantum Tori can be viewed as objects in this category, and we use this fact to calculate various fundamental groups and apply these results to compute morphisms between Hyper Quantum Tori.  Having obtained a notion of continuous maps between Hyper Quantum Tori, we use the preceeding work to define a notion of continuous maps between Quantum Tori, and prove the following result:

\begin{thm}\label{A}
Let $\mathcal{QT}$ denote the category of Quantum Tori, whose morphisms are continuous maps in the category of locally internal quotient spaces.  Then the categories $\mathcal{QT}$ and $\mathcal{PL}$ are equivalent.
\end{thm}

\S \ref{lift2} is an important discussion concerning the definition of morphisms in $\mathcal{LIQ}$ we gave in Definition \ref{morphisms}. In particular we note the limitations in the use of nonstandard mathematics in this area, and the possibility of using more sophisticated model theoretic techniques to gain better insight in to this subject.\\

The remainder of the paper concentrates on the arithmetic consequence of our definitions, and the potential use of Quantum Tori in solving Hilbert's twelfth problem for real quadratic fields.  An immediate consequence of Theorem \ref{A} is the following:

\begin{thm}\label{B}
Let $Z_L$ be a Quantum Torus.  Then the endomorphism ring of $Z_L$ is isomorphic to either $\B{Z}$ or an order of a real quadratic field.
\end{thm}

This result is analogous to a result for elliptic curves in the theory of Complex Multiplication, which states that the endormorphism ring of such an object is isomorhic to either $\B{Z}$ or an order in an \emph{imaginary} quadratic field.\\

We now fix a real quadratic field, $F$.  We may consider those Quantum Tori whose endomorphism ring is isomophic to the maximal order in $F$. This maximal order is realised as the ring of integers of $F$, and denoted $\mathcal{O}_F$.  Given our definiton of morphisms in $\mathcal{QT}$ it is easy to see that the following is well defined:

\begin{defn}
Let $F$ be a real quadratic field.  Define
$$\mathcal{QT}(\mathcal{O}_F) :=\frac{\textrm{Quantum Tori $Z$ with End($Z) \cong \mathcal{O}_F$}}{
\textrm{Isomorphism}}.$$
\end{defn}

We then prove the following result:

\begin{thm}\label{C}
Let $Z_L$ be a Quantum Torus, whose endomorphism ring is isomorphic to an order of a real quadratic field $F$.  Let $H_F$ denote the Hilbert class field of $F$, and $\textrm{Gal}(H_F/F)$ the galois group of the extension $H_F/F$.  Then there exists a simply transitive action of $\textrm{Gal}(H_F/F)$ on $\mathcal{QT}(\mathcal{O}_F)$.
\end{thm}

This last results represents an important step in the use of Quantum Tori in Real Multiplication, and we discuss the implications of this result.

\section{Nonstandard Models and Hyper Quantum Tori}\label{modeltheory}

Let $\mathcal{L}$ be the language $\{ \in\}$, and let $S$ be an infinite set.  We define the \emph{superstructure over $S$} by  
$$V(S) :=\bigcup_{n=0}^\infty V_n(S)$$
where $V_0(S):=S$ and $V_{n+1}(S)=V_n(S) \cup \mathcal{P}(V_n(S))$ for all $n \geq 1$.  Then $V(S)$ is an $\mathcal{L}$-structure.

\begin{defn}
Let ${}^\ast: V(S) \rightarrow V(R)$ be a map of superstructures.  We say that $V(R)$ is a \emph{nonstandard model} of $V(S)$ if the following properties are satisfied:
\begin{itemize}
\item ${^*S}=R$.
\item ${^*s}=s$ for all $s \in S$.
\item {\bf Transfer principal:} For every formula $\sigma(x_1, \ldots, x_n)$ with $n$-variables elements $a_1, \ldots, a_n \in V(S)$ we have
$$V(S) \models \sigma(a_1,\ldots, a_n) \Leftrightarrow V(R) \models \sigma({^*a_1},\ldots, {^*a_n}). $$
\item {\bf Saturation:} If $A \in V(S)$ is a collection of sets with the finite intersection property, then 
$$\bigcap \{{^*B}: B \in A\} \neq \emptyset. $$
\end{itemize}
\end{defn}

\begin{defn}[Internal sets] 
A set $A \in V({^*S})$ is called \emph{internal} if there exists $B \in V(S)$ such that $A \in {^*B}$.
\end{defn}
Equivalently, the internal sets of $V({^*S})$ are seen to be the definable sets in the $\mathcal{L}$-structure:
\begin{prop}
A set $A \in V({^*S})$ is internal if and only if there exists an $\mathcal{L}$-formula $\phi(\bar{x},\bar{y})$ and elements $a_1,\ldots, a_n$ such that
$$A=\{x: V({^*S}) \models \phi(x,\bar{a})\}.$$
\end{prop}
As a consequence of this we obtain property of nonstandard models known as \emph{permenance}:
\begin{prop}[Permenance]\label{permenance}
Let $\phi$ be an internal statement in $n$ variables.  Then the set of $\bar{a} \in {^*M^n}$ for which
${^*\mathcal{M}} \models \phi(\bar{a})$ is an internal set.
\end{prop}

We now list some terminology specific to studying topological space within a nonstandard structure.  We restrict our study to the case when $S=\B{R}$. 

\begin{defn}
Let $X \in V({\B{R}})$ be a topological space, whose topology consists of a collection of open sets $\mathcal{T}$.  Then we make the following definitions:
\begin{itemize}
\item An element of ${^*X}$ which lies in the image of ${^*}: X \rightarrow {^*X}$ is called \emph{standard}.
\item Let $x$ be an element of ${^*X}$, and suppose that $y \in {^*X}$ is such that $y \in {^*U}$ for all $U \in \mathcal{T}$. Then we say $y$ is \emph{infinitesimally close} to $x$, and write $x \simeq y$.
\item If $y \in {^*X}$ is infinitesimally close to a standard element we say $y$ is \emph{near standard}.
\end{itemize}
\end{defn}

\begin{defn}
If $X$ is a topological space in $\mathcal{M}$ we have two natural topologies on ${^*X}$:
\begin{itemize}
\item The $Q$-topology: If $\mathcal{B}$ is a basis of open sets for the topology on $X$, then ${^*\mathcal{B}}$ is a basis of open sets for the $Q$-topology on $^*X$;
\item The $S$-topology: A basis of open sets on ${^*X}$ is given by $\{ {^*A}: A \in \mathcal{B}\}.$
\end{itemize}
\end{defn}

The standard part map provides the bridge between the Nonstandard and ``standard'' models:

\begin{defn}[Standard part]
Let $X$ be a topological space, and let ${^*X^{\textrm{nst}}}$ denote the near standard elements of ${^*X}$.  Let $$\textrm{st}:{^*X^{\textrm{nst}}} \rightarrow X$$
denote the map sending $x \in {^*X^{\textrm{nst}}}$ to the unique standard element $x_0$ of $X$ such that $x \simeq x_0$.  When $X=\B{R}$ it can be shown that $\textrm{st}$ is an additive homomorphism, and if $x,y \in {^*\B{R}^{\textrm{nst}}}$ then $\textrm{st}(xy)=\textrm{st}(x)\textrm{st}(y)$.
\end{defn}

\subsection{Hyper Quantum Tori}\label{HQT}

Using the ideas of the previous section, we make the following definition:

\begin{defn}[Hyper Quantum Torus]
For a pseudolattice $L \subseteq \B{R}$, define the Hyper Quantum Torus associated to $L$ by $T_L:=\R/L$.  We let $\pi_L$ denote the natural projection $\R \rightarrow T_L$.
\end{defn}

\begin{prop}\label{haus}
Let $T_L$ be a Hyper Quantum Torus endowed with the induced quotient topology from the Q-topology on ${^*\B{R}}$.  Then $T$ is a Hausdorff space.
\end{prop}
\begin{proof}
We first observe that $L$ is not dense in $\R$.  Given $x \in \R$ consider the monad $\mu(x)=\{y \in \R: y \simeq x\}$.
Suppose $z \in (x+L )\cap \mu(x)$.  Then $z=x+l$ for some $l \in L$, but $x \simeq x+l$, and hence $l \simeq 0$.  Since every element of $L$ is standard we have $l=0$.  \\

Now suppose $x+L,\ y+L \in T_L$.  Because $\R$ is Hausdorff and $L$ is not dense in $\R$, there exist open sets $U_x$ and $U_y$ of $\R$ containing $x$ and $y$ respectively such that
$$
\begin{array}{c}
U_x \cap (x+L)=\{x\}\\
U_y \cap (y+L)=\{y\}\\
U_x \cap U_y = \emptyset.
\end{array}
$$
Now put $\mathcal{U}_x=U_x+L$ and $\mathcal{U}_y=U_y+L$.  These are open disjoint subsets of $T_L$ containing $x+L$ and $y+L$ respectively.
\end{proof}

\section{Locally Internal Topological Spaces}\label{monadicspace}

The proof of Proposition \ref{haus} shows that any point of $T_L$ has a neighbourhood which is isomorphic to an internal subset of $\R$ - if $x +L \in T_L$, then the set
$$\mathcal{V}_x^\varepsilon:=\{y+L: (y \in \R) \land (\left| y-x\right|<\varepsilon)\}$$
is isomorphic to the open interval $(-\varepsilon,\varepsilon)$ for any infinitesimal $\varepsilon$.  This inspires the following definition:

\begin{defn}[Locally Internal Topological Space]\label{lis}
Let $S$ be a topological space in a Nonstandard structure such that
\begin{quote}
For every $s \in S$ there exists an open neighbourhood $V_s$ of $S$ such that $V_s$ is isomorphic to 
an internal topological space.
\end{quote}
We say that $S$ is a \emph{locally internal} topological space.
\end{defn}

\subsection{Internal Covering Spaces}\label{covermap}

Suppose $X$ is a standard topological space.  Basic results in topology \cite{Armstrong} imply that if $\gamma$ is a path in $X$, and $\tilde X$ is a covering space for $X$, then $\gamma$ lifts to a path $\tilde \gamma$ in $\tilde X$.  This result implies that a continuous map between topological spaces $X$ and $Y$ lifts to a map between covers $\tilde X$ and $\tilde Y$ of these respective spaces.\\

We wish to have an analogous situation for locally internal topological spaces, where the covering space is an \emph{internal} topological space.  

\begin{defn}[Internal Covering space]\label{cover}
An internal cover of a locally internal space $S$ is a pair $(\tilde S,p)$ such that
\begin{itemize}
\item $\tilde S$ is an internal topological space;
\item $p$ is a surjective map from $\tilde S$ to $S$ satisfying the following condition:
\begin{quote}
For every $s \in S$, there exists an open neighbourhood $U_s$ of $s$; an isomorphism $\psi_s$ of $U_s$ on to an internal topological space $\mathfrak{U}_s$, and a decomposition of $p^{-1}(U_s)$ as a family $\{V_{s,i}\}$ of disjoint open internal subsets of $\tilde S$ such that the restriction of $\phi_s \circ p$ to $V_{s,i}$ is an internal homeomorphism from $V_{s,i}$ to $\mathfrak{U}_s$.
\end{quote}
\end{itemize}
We say that $S$ is a \emph{locally internal quotient} space.
\end{defn}

Internal covering spaces are not unique, as the following example shows:

\begin{quote}
Let $L=\B{Z}\omega_1+\B{Z}\omega_2$ be a pseudolattice, and let ${^*S^1} \cong {\R/\B{Z}\omega_1}$ denote the unit circle.  Consider the pair $({^*S^1},p_1)$, where $p_1$ is the map
$$
\begin{array}{rcl}
p_1:{^*S^1} & \rightarrow & T_L\\
x+\B{Z}\omega_1 & \mapsto & x+L.
\end{array}
$$
Then $({^*S^1},p_1)$ is an internal covering space for $T_L$.\\

Now let $\mathfrak{L}:{^*\B{R}} \rightarrow \R$ be a map of the form
$$\mathfrak{L}(x)=x+m(x)\omega_1$$
where $m(x)$ is some integer (depending on $x$) such that $x+m(x)\omega_1$ is limited.  Then the pair $({^*\B{R}},p_2)$ is an internal covering space for $T_L$, where
 $$
\begin{array}{rcl}
p_2:{^*\B{R}} & \rightarrow & T_L\\
x & \mapsto & \mathfrak{L}(x)+L.
\end{array}
$$
\end{quote}

\subsection{Morphisms between locally internal quotient spaces}\label{morphhqt}

Given a Hausdorff locally internal space $S$, in general $S$ will be an external object in our Nonstandard structure.  As a consequence the space of continuous functions on such a space can be very big.  When $S$ has an internal cover, we use this to restrict the space of such functions by the following definition:

\begin{defn}\label{morphisms}
Let $S$ and $T$ be locally internal quotient spaces. Let $p_S : \tilde S \rightarrow S$
and $p_T : \tilde T \rightarrow T$ denote the respective covering maps.  A morphism between
$S$ and $T$ is a map $f:S \rightarrow T$ such that there exists an internal function
$\tilde f : \tilde S \rightarrow \tilde T$ such that 
$$p_T \circ \tilde f(\tilde s) = f \circ p_S(\tilde s)$$
for all $\tilde s \in \tilde S$.  We say that a morphism $f$ is Q-continuous if
$\tilde f$ is Q-continuous.

\end{defn}
From this definition it is clear that the composition of two morphisms is again a morphism.\newline
\begin{defn}[Category of Locally Internal Quotient spaces]
Let $\mathcal{LIQ}$ be the category such that 
\begin{itemize}
\item The objects of $\mathcal{LIQ}$ are locally internal quotient spaces;
\item A morphism between locally internal quotient spaces $S$ and $T$ is as defined in Definition \ref{morphisms}.
\end{itemize}
\end{defn}

\subsection{The fundamental group of a locally internal space}\label{fundom}

For a standard topological space $X$, a path in $X$ is a continuous map $\gamma: I \rightarrow X$ where $I$ is the unit interval.  Taking the $\ast$-transform of this definition, an internal path in an internal topological space $Y$ in a Nonstandard structure is a Q-continuous map $\gamma:{^*I} \rightarrow Y$ where
$${^*I}:=\{x \in {^*\B{R}} : 0 \leq x \leq 1 \}.$$
We note that ${^*I}$ is trivially a locally internal quotient space (covered by itself) and extend the above notion to define paths in objects in $\mathcal{LIQ}$.

\begin{defn}
Let $S$ be a locally internal topological space.
A path in $S$ is a Q-continuous morphism $\gamma$ in $\mathcal{LIQ}$ from ${^*I}$ to $S$.
We say that a path $\gamma$ is a loop based at $s \in S$ if $\gamma(0)=\gamma(1)=s$.
\end{defn}

Similarly we can extend the notion of homotopies between paths in locally internal quotient spaces:

\begin{defn}
Let $\gamma_1$ and $\gamma_2$ be paths in a locally internal quotient space $S$.  A homotopy $F$
between $\gamma_1$ and $\gamma_2$ is a Q-continuous morphism in $\mathcal{LIQ}$ from ${^*I^2}$ to $S$
such that
\begin{itemize}
\item $F(t,0)= \gamma_1(t)$ for all $t \in {^*I}$;
\item $F(t,1)= \gamma_2(t)$ for all $t \in {^*I}$.
\end{itemize}
We say the two paths $\gamma_1$ and $\gamma_2$ are homotopic and write $\gamma_1 \simeq \gamma_2$.
If we wish to refer explicitly to the homotopy $F$ we may write $\gamma_1 \simeq_F \gamma_2$.
If $\gamma_1$ and $\gamma_2$ agree on some subset $A$ of ${^*I}$, we say that $F$ is a homotopy between
$\gamma_1$ and $\gamma_2$ relative to $A$ if we have the additional condition
\begin{itemize}
\item $F(a,s)= \gamma_1(a)=\gamma_2(a)$ for all $a \in A$.
\end{itemize}
\end{defn}

It is easily shown that the relation $\simeq$ is an equivalence relation. If $\gamma$
is a path in $S$ we let $\langle \gamma \rangle$ denote the equivalence (or homotopy) class of $\gamma$.
Given $s \in S$, we denote the set of homotopy classes relative to $\{0,1\}$ of loops in $S$
based at $s$ by $\pi_1(S,s)$.
We define a law of composition on $\pi_1(S,s)$ by
$$\langle \gamma_1 \rangle.\langle \gamma_2 \rangle = \langle \gamma_1 \star \gamma_2 \rangle$$
where
$$\gamma_1 \star \gamma_2(t) = \left\{
\begin{array} {ll}
\gamma_1(2t) & \textrm{ if $0 \leq t \leq \frac{1}{2}$}\\
\gamma_2(2t-1) & \textrm{ if $\frac{1}{2} \leq t \leq 1$.}
\end{array} \right.
$$
One can prove using exactly the same methods as for standard topological spaces that $\pi_1(S,s)$ is a group under this operation \cite{Armstrong}. The identity element is the constant loop at $s$. \newline

\noindent
\begin{rmk} \label{rtop} Observe that if $S$ is an internal topological space then this agrees with
the natural definition (the $\ast$-transform of the standard definition) of $\pi_1(S)$.
At first it may seem to be a stronger definition since we have the property that paths
and homotopies lift to internal covers, but these results follow for internal spaces by
$\ast$-transform of the standard results. 
\end{rmk}

Suppose $f: S \rightarrow T$ is a morphism of locally internal quotient spaces.  Then $f$ induces a homomorphism
$f_\ast: \pi_1(S,s)  \longrightarrow  \pi_1(T,f(s))$.\\

\begin{prop}
Let $(\tilde S,p)$ be an internal covering space for a locally internal space $S$.  If $\tilde S$ is path connected then for any $\tilde s \in \tilde S$ the induced map
$p_\ast: \pi_1(\tilde S,\tilde s) \rightarrow \pi_1(S,s)$ is injective, where $s=p(\tilde s)$.
\end{prop}
\begin{proof}
Suppose $\tilde \gamma$ is a loop in $\tilde S$ such that $\gamma:=p \circ \tilde \gamma$
is null homotopic.  Choose a specific homotopy $1_s \simeq_F \gamma$, where $1_s$ denotes
the constant loop at $s$.  Choose a Q-continuous lift $\tilde F$ of $F$ such that
$p \circ \tilde F=F$.  We may assume that $\tilde F(0)=\tilde s$ since if not, chose a path $p$ from $\tilde F(0)$ to $\tilde s$ and replace $\tilde F$ by the homotopy $G$ such that for each $s \in {^*I}$
$$\tilde G(t, s)=p^{-1} \star \tilde F( \ ,s) \star p(t).$$  \newline

Note that once we have fixed $\tilde F(0)$ the function $\tilde F$ is unique.  Suppose there were two such lifts $\tilde F_1$ and $\tilde F_2$. Then we would have $\tilde F_1(t)-\tilde F_2(t) \in \ker(p)$.  The kernel of $p$ is a Q-discrete set since for each point $s$, $p$ identifies an internal neighbourhood of $s$ homeomorphically with an internal subset of $\tilde S$.
Since both the lifts are Q-continuous and internal maps of the connected
set ${^*I}$ and agree at $0$, we have $\tilde F_1=\tilde F_2$. \newline

We need to show that $\tilde F$ gives a homotopy from $\tilde \gamma$ to the constant
loop at $\tilde s$.  Let $P$ denote the internal path connected set
$$\{(t,0) \in {^*I}^2 : 0 \leq t \leq 1\} \cup \{(0,t) \in {^*I}^2 : 0 \leq t \leq 1\}.$$
Since $F$ is a homotopy relative to $\{0,1\}$ we see that $F(P)=s$. Since
$p \circ \tilde F = F$ we have $\tilde F(P) \in p^{-1}(s)$, which as we have seen is
a Q-discrete set.  Since $\tilde F$ is internal and Q-continuous we have $\tilde F(P)=\tilde s$.
This shows that the path $F(t,0)$ is the constant loop at $\tilde s$.  The path
$F(t,1)$ is a lift of $\gamma$ which starts at $\tilde s$.  There is a unique such path by
an analogous argument to the uniqueness of homotopy lifting in the above paragraph.  Hence
$F(t,1)=\tilde \gamma(t)$.

\end{proof}

We stress that despite working within a nonstandard model, the proof of the previous proposition does not require any new ingredients mathematically.  Once we have the properties of path and homotopy lifting we are working with internal functions on the covering space and the proofs carry through by $\ast$-transform of the standard results. 

\begin{prop}    \label{maplifting}
Let $f: S \longrightarrow T$ be a morphism between monadically internal spaces, let $s \in S$ and suppose that $\tilde S$ and $\tilde T$ are path connected.
There is a lift $\tilde f: S \longrightarrow \tilde T$ such that $f(s)=\tilde t$ if and only if
$f_\ast(\pi_1(S,s)) \subseteq \pi_\ast(\pi_1(\tilde T,\tilde t))$. This lift is unique.
\end{prop}
\begin{proof}
The result is proved for standard spaces in \cite{Armstrong}. It is easy to adapt these techniques to obtain the result for locally internal quotient spaces.
\end{proof}

\begin{defn}[Covering Transformation]
A covering transformation for an internal covering space $(\tilde S,\pi)$ is an internal homeomorphism $h: \tilde S \rightarrow \tilde S$ such that $\pi \circ h = \pi$.
\end{defn}

The set of all covering transformations forms a group $\textrm{Cov}(\tilde S/S)$
under composition, and we have an action of $\textrm{Cov}(\tilde S/S)$ on $\tilde S$ by 
$$
\begin{array}{rcl}
\textrm{Cov}(\tilde S/S) \times \tilde S & \rightarrow & \tilde S \\
(h,\tilde s) & \mapsto & h(\tilde s).
\end{array}
$$
Note that
if $h_1$ and $h_2$ are covering transformations which agree on a point $\tilde s$, then
both $h_1(x)-h_2(x)$ and the constant map $K(x)=0$ take the value $0$ at $x=\tilde s$ and lift
$\pi$.  Hence by the uniqueness part of Proposition \ref{maplifting} we have $h_1=h_2$.\\

Utilising methods from the proof of the corresponding standard result, it is easy to show that the following is true:

\begin{prop} \label{covertrans}
Let $(\tilde S, p)$ be a path connected internal covering space.
Suppose $p_\ast(\pi_1(\tilde S,\tilde s))$ is a normal subgroup of $\pi_1(S,s)$.  Then
$\textrm{Cov}(\tilde S/S)$ is isomorphic to the quotient $\pi_1(S,s) /p_\ast(\pi_1(\tilde S,\tilde s))$.
\end{prop}

We now use this result to describe various fundamental groups associated to Hyper Quantum Tori.

\subsection{Covering spaces for Hyper Quantum Tori}\label{coverqt2}

The definitions of the previous sections enable us to calculate the fundamental group associated to Hyper Quantum Tori:

\begin{prop} \label{iso}
Let $L$ be a pseudolattice.  For any $z \in T_L$, we have an isomorphism
$\pi_1(T_L,z) \cong {^*\B{Z}} \oplus \B{Z}$.
\end{prop}
\begin{proof}
Consider the internal covering space $({^*\B{R}},p_2)$ considered in \S\ref{covermap} :
$$
\begin{array}{rcl}
p_2:{^*\B{R}} & \rightarrow & T_L\\
x & \mapsto & \mathfrak{L}(x)+L
\end{array}
$$
Since $\pi_1({^*\B{R},r})$ is trivial for any $r \in {^*\B{R}}$ (see Remark \ref{rtop}), Proposition
\ref{covertrans} gives an isomorphism
$\textrm{Cov}({^*\B{R}}/T_L) \cong \pi_1(T_L,z)$.
Suppose $L=\B{Z}\omega_1+\B{Z}\omega_2$.  For each $l\in {^*\B{Z}}\omega_1 + \B{Z}\omega_2$
let $h_{l}$ be the covering transformation defined by $h_{l}(r)=r+l$, and consider the homomorphism
$$
\begin{array}{rcl}
{^*\B{Z}}\omega_1 + \B{Z}\omega_2 & \rightarrow & \textrm{Cov}({^*\B{R}}/\mathfrak{F}_\theta) \\
l & \mapsto & h_{l}.
\end{array}
$$
By the previous discussion we see that elements $f \in \textrm{Cov}({^*\B{R}}/T_L)$
are determined by their value at $0$, hence this map is injective.  Given
$f \in  \textrm{Cov}({^*\B{R}}/T_L)$ we have $f(0) \in {^*\B{Z}}\omega_1+\B{Z}\omega_2$,
and hence $f=h_{f(0)}$.  Hence the above map defines an isomorphism ${^*\B{Z}}\omega_1 + \B{Z}\omega_2 \cong \textrm{Cov}({^*\B{R}}/\mathfrak{F}_\theta)$.  Finally note that ${^*\B{Z}}\omega_1+\B{Z}\omega_2 \cong {^*\B{Z}} \oplus \B{Z}$.
\end{proof}

When we come to consider how morphisms between Hyper Quantum Tori give us maps between Quantum Tori, we will want to consider the standard part images of covering transormations.  However, by Proposition \ref{iso} such transformations are translations by elements of ${^*\B{Z}\omega_1}+\B{Z} \omega_2$, and as such there exist some covering transformations which do not map $\R$ to itself.  We define the following subcategory of $\mathcal{LIQ}$:

\begin{defn}
Let $\mathcal{LIQ}^{\textrm{nst}}$ be the subcategory of $\mathcal{LIQ}$ such that 
\begin{itemize}
\item The objects of $\mathcal{LIQ}^{\textrm{nst}}$ are objects of $\mathcal{LIQ}$ such that 
\begin{quote}
The restriction of the projection $p: \tilde S \rightarrow S$ to those limited elements of $\tilde S$ is an S-continuous surjection on to $S$. 
\end{quote}
\item  Morphisms in $\mathcal{LIQ}^{\textrm{nst}}$ are those morphisms in $\mathcal{LIQ}$ which map limited elements to limited elements.
\end{itemize}
\end{defn}
For $z \in T_L$, we let $\pi_1^{\textrm{nst}}(T_L,z)$ denote the fundamental group of $T_L$ based at $z$ where all the paths and homotopies are required to be morphisms in $\mathcal{LIQ}^{\textrm{nst}}$.
\begin{prop} \label{limfun}
For any $z \in T_L$, we have an isomorphism $\pi_1^{\textrm{nst}}(T_L,z) \cong L$.
\end{prop}
\begin{proof}
Proposition \ref{iso} implies that the following injection is an isomorphism:
$$
\begin{array}{rcl}
\Phi:\textrm{Cov}({^*\B{R}}/T_L) & \rightarrow & \pi_1(T_L,z)\\
h & \mapsto & \gamma_h
\end{array}
$$
where
$$\gamma_h(t):=\pi \circ (\tilde z+t(h(1)-h(0)). $$
Such a transformation preserves $\R$ if and only if it is a translation by an element of $\B{Z}\omega_1+\B{Z}\omega_2=L$. 
\end{proof}

Hence we may obtain the pseudolattice $L$ from the fundamental group of the Hyper Quantum Torus $T_L$.  We view this as analogous to the determination of the pseudolattice $L_\theta=\B{Z}+\theta \B{Z}$ from the K-theory of the Noncommutative Torus $A_\theta$ in Noncommutative Geometry \cite{Manin}.

\subsection{Defining morphisms between Quantum Tori}   \label{endomorphism}

Our motivation for defining and studying locally internal quotient spaces arose from the observation that a Hyper Quantum Torus $T_L$ can be viewed as an object in the category of such spaces.  We have developed a notion of morphism in between such spaces, which we now apply specifically to Hyper Quantum Tori:

\begin{defn}[Category of Hyper Quantum Tori]
Let $\mathcal{HQT}$ denote the category such that 
\begin{itemize}
\item The objects of $\mathcal{HQT}$ are Hyper Quantum Tori;
\item The morphisms between Hyper Quantum Tori correspond to those homomorphisms of the universal internal covering space of Hyper Quantum Tori which map limited elements to limited elements. 
\end{itemize}
\end{defn}

\begin{lemma}\label{mo}
We have a bijection between the set of morphisms $T_L \rightarrow T_M$ in $\mathcal{HQT}$, and the set of nonzero real numbers $\alpha$ such that $\alpha(L) \subseteq M$.
\end{lemma}
\begin{proof}
Proposition \ref{maplifting} shows that a covering space for $T_L$ is universal if it has trivial fundamental group.  Hence ${^*\B{R}}$ is the universal internal covering space for Hyper Quantum Tori.  Definition \ref{morphisms} implies that the morphisms between Hyper Quantum Tori are homomorphisms $\phi: {^*\B{R}} \rightarrow {^*\B{R}}$ which map $\R$ to itself and satisfy the following condition:
\begin{equation}\label{proj}
\phi(r)+M=\phi(r+L)\qquad  \forall r \in \R.
\end{equation}
Since $\phi$ is an internal homomorphism of ${^*\B{R}}$ it is equal to multiplication by $\alpha$ for some $\alpha \in {^*\B{R}^\times}$.  By \eqref{proj} we have $\alpha(L) \subseteq M$, which since $L$ and $M$ are standard imply that $\alpha \in \B{R}^\ast$.  Conversely, multiplication by any such element induces a morphism in $\mathcal{HQT}$.

\end{proof}

\begin{cor}
The categories $\mathcal{PL}$ and $\mathcal{HQT}$ are equivalent.
\end{cor}
\begin{proof}
This follows from the above result, and the definition of $\mathcal{PL}$ in Definition \ref{PL}.
\end{proof}

The natural projection $\pi_L: \R \rightarrow T_L$ induces a well defined map $p_L$ on $T_L$ whose image is a Quantum Torus:
$$
\begin{array}{rcl}
p_L: T_L & \rightarrow & Z_L\\
\pi_L(x) & \mapsto & \textrm{st}(x)+L. 
\end{array}
$$

\begin{defn}\label{moqt}
A morphism between Quantum Tori $Z_{L}$ and $Z_{M}$ is a map $f:Z_{L} \rightarrow Z_{M}$ such that there exists a map $\tilde f: T_{L} \rightarrow T_{M}$ such that the following diagram commutes:
$$
\xymatrix{
T_{L} \ar[d]^{p_L} \ar[r]^{\tilde f} & T_M \ar[d]^{p_M} \\
Z_L \ar[r]^{f} & Z_M.
}
$$\end{defn}
\begin{proof}[Proof of Theorem \ref{A}]
Since every morphism of Hyper Quantum Tori is S-continuous, every element morphism of Hyper Quantum Tori yields a morphism of Quantum Tori.  We need to show that this association is injective.\\

Note that for each pseudolattice we have an injection $\iota_L: Z_L \hookrightarrow T_L$, which satisfies $p_L \circ \iota_L=1_{Z_L}$.  We need to prove that if $g:T_L \rightarrow T_M$ is a morhpism of Hyper Quantum Tori such that
\begin{equation}\label{asso}
p_M \circ g \circ \iota_L(x)=x
\end{equation}
for all $x \in Z_L$, then $g$ is the identity.\\

Since the induced map on Quantum Tori is the identity we have $L=M$.  By Proposition \ref{mo}, $g$ is equal to multiplication by $\alpha$ for some $\alpha \in \B{R}^\times$ such that $\alpha L \subseteq L$.  Then \eqref{asso} implies that
$$\alpha r+L=r+L$$
for all $r \in \B{R}$, and hence $\alpha=1$.
\end{proof}
\subsection{Limitations of Nonstandard Analysis}\label{lift2}

In Definition \ref{morphisms} we defined morphisms between locally internal quotient spaces to possess the property that they lifted to internal maps between their internal covering spaces.  In this section we look at whether it is possible to obtain an equivalent definition without reference to such a cover, which would allow us to consider maps between the more general class of locally internal topological spaces.  Rather than consider the general case of maps between locally internal quotient spaces, we shall consider the specific problem of defining paths in such spaces. 

\begin{defn}\label{loc}
Let $S$ be a locally internal space.  A \emph{locally internal} path in $S$ is a Q-continuous map $\gamma : {^*I} \rightarrow S$
such that for all $t \in {^*I}$, there exists an internal neighbourhood $U_t$ of $t$ with the following properties:
\begin{enumerate}
\item The image of $U_t$ lies in an open neighbourhood $V_{\gamma(t)}$ of $\gamma(t)$;
\item There exists an isomorphism $\phi_{\gamma(t)}: V_{\gamma(t)} \cong \mathfrak{V}_{\gamma(t)}$ for some internal topological space $\mathfrak{V}_{\gamma(t)}$;
\item The composition $\phi_{\gamma(t)} \circ \gamma: U_t \rightarrow \mathfrak{V}_{\gamma(t)}$ is an internal map. 
\end{enumerate}
\end{defn}

Hyper Quantum Tori posses a slightly stronger property than local internality.  Instead of the condition of Definition \ref{lis} we have 
\begin{quote}
For every $z \in T_L$, for every infinitesimal $\varepsilon$ the set $(z-\varepsilon,z+\varepsilon)+L \in T_L$ is an open neighbourhood of $z$ isomorphic to an internal topological space.
\end{quote}
Consequently we shall work with a stronger notion of paths defined by the following:
\begin{defn}\label{mon}
A \emph{monadically internal} path in $T_L$ is a map $\gamma: {^*I} \rightarrow T_L$ such that (with the notation of Definition \ref{loc}), for every $t \in {^*I}$, and every infinitesimal $\varepsilon$, there exist $V_{\gamma(t)}$, $\mathfrak{V}_{\gamma(t)}$, $\phi_{\gamma(t)}$ corresponding to $U_{t}=(t-\varepsilon,t+\varepsilon)$.
\end{defn}

When $\gamma$ is a monadically internal path in $T_L$, we can deduce some information about a lift of $\gamma$ to ${^*\B{R}}$:

\begin{prop}    \label{Smonad}
Fix $t \in {^*I}$.  Let $\gamma$ be a monadically internal path in $T_L$, and let $\tilde x \in {^*\B{R}}$ be such that $\pi(\tilde x)=\gamma(t)$.  Then there exists a unique internal function $f$ and $r \in \B{R}$ such that $f: [t-r,t+r]\cap {^*I} \rightarrow {^*\B{R}}$, $f(0)=\tilde x$ and on $\mu(t)$ we have $\pi \circ f = \gamma$.
\end{prop}
\begin{proof}
This is a consequence of the property of \emph{permanence} introduced in Proposition \ref{permenance}. We first suppose that $t \notin \mu(0) \cup \mu(1)$.\newline

For each $\eta \simeq 0$, there exist internal sets $\mathfrak{V}_{\gamma(t)}^\eta \subseteq {^*\B{R}}$ and functions $\phi_{\gamma(t)}^\eta$ such that $\pi \circ \phi_{\gamma(t)}^\eta(t)=\gamma(t)$ and 
$$\phi_{\gamma(t)}^\eta: (t-\eta,t+\eta) \rightarrow \mathfrak{V}_{\gamma(t)}^\eta$$ is an internal function.\\

Fix $\eta_0 \simeq 0$.  Then for each $\eta \simeq 0$ with $\eta > \eta_0$ we obtain an internal function $\phi_{\gamma(t)}^\eta$ such that on $(t-\eta_0,t+\eta_0)$ we have $\phi_{\gamma(t)}^\eta=\phi_{\gamma(t)}^{\eta_0}$. By permanence there exists $r \in \B{R}$ and an internal function $\phi_{\gamma(t)}^r$ and an internal open set $\mathfrak{V}_{\gamma(t)}^r$ such that 
$\phi_{\gamma(t)}^{r}: (t-r,t+r) \rightarrow \mathfrak{V}_{\gamma(t)}^r$
and the restriction of $\phi_{\gamma(t)}^{r}$ to $(t-\eta_0,t+\eta_0)$ is $\phi_{\gamma(t)}^{\eta_0}$.
Since this holds for any $\eta_0 \simeq 0$ we see that $\phi_{\gamma(t)}^{r}$ agrees with $\phi_{\gamma(t)}^{\eta_0}$ for all $\eta_0 \simeq 0$.   \newline

Now suppose we had two such lifts $\phi_{\gamma(t)}^{r}$ and $\psi_{\gamma(t)}^{r}$ defined on $(t-r,t+r)$.  Consider the internal set
$$S:=\{t' \in (t-r,t+r) :\phi_{\gamma(t)}^{r}(t)=\psi_{\gamma(t)}^{r}(t)\}.$$
Consider the internal statement
$$\varphi(s): \ \ t' \in (t-s,t+s) \Rightarrow t' \in S.$$
Then $\psi(\varepsilon)$ is valid for all $\varepsilon \simeq 0$, hence by permanence there exists
$r' \in {\B{R}}$ such that $\psi(r')$ holds.
Hence $\phi_{\gamma(t)}^{r}$ is unique on $(t-r',t+r')$.\\

If $t \in \mu(0)$ then we apply the above techniques to the interval $[0,\eta)$, and similarly if $t \in \mu(1)$ we consider the interval $(\eta,1].$
\end{proof}

Note that the above proof does \emph{not} give us a lift of $\gamma$ on the whole interval $(t-r,t+r)$.  If $T_L$ is a Hyper Quantum Torus, then the path given by
$$
\begin{array}{rcl}
\gamma: {^*I} & \rightarrow & T_L\\
t & \mapsto & \textrm{st}(t)+L
\end{array}
$$
defines a monadically internal path in $T_L$, but which has no lift to an internal map $f: {^*I} \rightarrow {^*\B{R}}$.
This demonstrates that Definition \ref{loc} does not provide a notion of paths in locally internal spaces which lift to internal paths in internal covering spaces. \\

By Proposition \ref{Smonad}, an equivalent definiton of a monadically internal path is given by:

\begin{defn}
A monadically internal path in $T_L$ is a function $\gamma:{^*I} \rightarrow F_L$
such that
\begin{quote}
For each $x \in {^*I} \subset \R$, for every $\varepsilon \simeq 0$ there exists
a Q-continuous internal function $\tilde \gamma_x: [x-\varepsilon,x+\varepsilon] \cap {^*I} \rightarrow {^*\B{R}}$ such that
for all $t \in  [x-\varepsilon,x+\varepsilon]\cap {^*I} $, $\pi \circ \tilde \gamma(t) =P^{-1} \circ  \gamma(t)$.
\end{quote}
\end{defn}

We know that this is not enough to give us the property of path lifting.  A natural weakening of this notion is given by 

\begin{defn}\label{app}
An appreciably internal path in $T_L$ is a function $\gamma:{^*I} \rightarrow F_L$
such that
\begin{quote}
For each $x \in {^*I} \subset \R$, there exists $r_x \in \B{R}$ and a Q-continuous internal function $\tilde \gamma_x: [x-r_x,x+r_x]\cap {^*I}  \rightarrow {^*\B{R}}$ such that
for all $t \in  [x-\varepsilon,x+\varepsilon]\cap {^*I} $, $\pi \circ \tilde \gamma(t) =P^{-1} \circ  \gamma(t)$.
\end{quote}
\end{defn}

This is a failure to define paths in $T_L$ with the path lifting property without reference to the covering space ${^*\B{R}}$.  However, it is easy to see that the compactness of $I$ implies that appreciably internal paths yield have internal lifts in internal covering spaces.\\

This highlights the limitations in the use of Nostandard Analysis in studying Quantum Tori - locally definable maps  between Hyper Quantum Tori do not lift to definable maps between covering spaces.   It would be interesting to study structures with fewer definable sets in which Quantum Tori can be interpreted to see whether it is possible to find a structure inwhich locally definable maps lift to definable maps of covers.

\section{Quantum Tori with Real Multiplication} \label{RM}

The theory of Complex Multiplication that forms the basis for the solution of Kronecker's Jugendtraum relies on the existence of elliptic curves over $\B{C}$ whose endomorphism ring is strictly greater than $\B{Z}$.  We let $\mathcal{QT}$ denote the category whose objects are Quantum Tori, and whose morphisms are continuous homomorphisms as described in Definition \ref{moqt}.  Together with the work of the previous section, we are now in a position to prove Theorem \ref{B}: \newline

\begin{proof}[Proof of Theorem \ref{B}]
Let $Z_L$ be a Quantum Torus.  By Theorem \ref{A}, the endomorphism ring of $Z_L$ is isomorphic to the set of those $\alpha \in \B{R}^\ast$ such that $\alpha L \subseteq L$.  Suppose there exists such an $\alpha$ such that $\alpha \notin \B{Z}$, and let $L=\B{Z}\omega_1+\B{Z}\omega_2$.  Then there exist $a,b,c$ and $d \in \B{Z}$ such that
\begin{equation}  \label{eq}
\begin{array}{rcl}
\alpha \omega_1&=&a\omega_1+b \omega_2\\
\alpha \omega_2&=&c\omega_1+d \omega_2.
\end{array}
\end{equation}
Dividing the second of these equations by $\omega_2$, we observe that since $\alpha \notin \B{Z}$, $c \neq 0$. We observe that $\theta:=\omega_1/\omega_2$ satisfies the quadratic equation
$$cX^2+(d-a)X-b=0. $$
Hence $[\B{Q}(\theta):\B{Q}]=2$, and $\B{Q}(\theta)$ is a real quadratic field.\\

We therefore conclude that $\textrm{End}(Z_L)$ is an integral extension of $\B{Z}$.  Eliminating $\alpha$ from \eqref{eq} we see that $\alpha$ satisfies the equation
$$X^2-(a+d)X+ad-bc=0.$$
Hence $\alpha$ is integral over $\B{Z}$ and therefore contained in the ring of integers $\mathcal{O}_F$ of $F=\B{Q}(\theta)$.\newline

We can therefore identify the ring of endomorphisms as a subring of the ring of integers of $F$. Hence $\textrm{End}(Z_L)$ is finitely generated as a $\B{Z}$-module and satisfies $\textrm{End}(Z_L) \otimes \B{Q} \cong F$.  These are the precisely the requirements for $\textrm{End}(Z_L)$ to be an order in $F$.
\end{proof}

\begin{defn}
Let $Z_L$ be a Quantum Torus such that $\textrm{End}(Z_L)$ is isomorphic to an order in a real quadratic field $F$.  We say that $Z_L$ has \emph{Real Multiplication} (by $F$).  We sometimes abbreviate this to say that $Z_L$ has RM.
\end{defn}

In this section we consider those Quantum Tori $Z_L$ such that $\textrm{End}(Z_L)$ is isomorphic to the \emph{maximal} order in a predefined real quadratic field $F$.  From a number theoretic point of view such orders have a special significance, being the ring of integers of $F$.  
 
\begin{lemma}    \label{frac}
Let $L$ be a pseudolattice in $\B{R}$. Then $\textrm{End}(Z_L) \cong \mathcal{O}_F$ if and only if
$L$ is homothetic to a fractional ideal in $F$.
\end{lemma}
\begin{proof}
Suppose $\textrm{End}(Z_L) \cong \mathcal{O}_F$.   Then
$$\mathcal{O}_F \cong \{x \in \B{R}: x L \subseteq L\}.$$
Suppose $L=\B{Z}\omega_1+\B{Z}\omega_2$. Then $L=\omega_1L_\theta$, where $\theta=\omega_2/\omega_1$, so
$$\mathcal{O}_F \cong \{x \in \B{R}: x L_\theta \subseteq L_\theta\}.$$
By the proof of Theorem \ref{B} we deduce that $\theta \in F$, and that we can identify $\mathcal{O}_F$ precisely with the set
$\{x \in \B{R}: x L_\theta \subseteq L_\theta\}$.  Hence $L$ is homothetic to a rank two $\mathcal{O}_F$-module contained in $F$. The latter object is precisely the definition of a fractional ideal in $F$.    \newline

Conversely let $\mathfrak{a}$ be a fractional ideal of $F$.  This is a rank two abelian subgroup of $\B{R}$, and therefore
a pseudolattice, so we may consider the quantum torus $Z_\mathfrak{a}$.  Every element of $\textrm{End}(Z_\mathfrak{a})$ lifts to  multiplication by $\alpha_\mathfrak{a}$ on $\B{R}$ such that $\alpha_\mathfrak{a} \mathfrak{a} \subseteq \mathfrak{a}$ for some $\alpha_\mathfrak{a} \in \B{R}$.  So $\textrm{End}(Z_\mathfrak{a})$ contains $\mathcal{O}_F$ as a suborder.  We know that $\textrm{End}(Z_\mathfrak{a})$ is isomorphic to an order of $F$, and that $\mathcal{O}_F$ is the maximal order.  Hence $\textrm{End}(Z_\mathfrak{a}) \cong \mathcal{O}_F$.
\end{proof}

\subsection{The algebraic nature of Quantum Tori}

Let $\sigma$ be a generator of the Galois group of $F$ over $\B{Q}$, and denote by $I_F$ the group of fractional
ideals of $F$. We have a natural inclusion $\iota : F^\ast \rightarrow I_F$ which sends $x$ to the principal ideal $(x)$
generated by $x$.
The class group of $F$ is a finite group \cite{Neukirch} defined by the quotient
$$C(F):=\frac{I_F }{ \iota(F^\ast)}. $$
The equivalence class of $\mathfrak{a} \in I_F$ in $C(F)$ is denoted by $[\mathfrak{a}]$. \newline

Let $$F^+:=\{x \in F: (x>0) \land (x^{\sigma}>0)\}$$
 denote the subgroup of totally positive elements of $F$.
The narrow class group of $F$ is defined to be
$$C(F)^+:=\frac{I_F }{\iota(F^+)}. $$
The equivalence class of $\mathfrak{a} \in I_F$ in $C(F)^+$ is denoted by $[\mathfrak{a}]^+$.  There is
a canonical surjection
\begin{equation}\label{narrowsurjection}
C(F)^+ \longrightarrow C(F).
\end{equation}

\begin{thm}  \label{action}
There is a well defined action of $C(F)$ on $\mathcal{QT(O}_F)$. This action is simply
transitive.
\end{thm}
\begin{proof}
Let $Z_L$ be a Quantum Torus with endomorphism ring isomorphic to $\mathcal{O}_F$.  Let $\mathfrak{a}$ be a fractional
ideal of $F$, and define $\mathfrak{a} \ast Z_L:=Z_{\mathfrak{a}L}.$
I claim this induces a well defined action of $I_F$ on $\mathcal{QT(O}_F)$.
\begin{itemize}
\item $\mathfrak{a}L$ is a pseudolattice. By Lemma \ref{frac} $L=\lambda \mathfrak{c}$ for some
$\lambda \in \B{R}^\ast$ and fractional ideal $\mathfrak{c}$ of $F$, so $\mathfrak{a}L=\lambda \mathfrak{ac}$.
The fractional ideal $\mathfrak{ac}$ is a pseudolattice since it is a rank two abelian subgroup of $\B{R}$, and hence
$\mathfrak{a}L$ is.
\item $\textrm{End}(Z_{\mathfrak{a}L}) \cong \mathcal{O}_F$.
With the notation of the last paragraph, $\mathfrak{a}L$ is homothetic to the fractional ideal $\mathfrak{ac}$.
The statement follows from Lemma \ref{frac}.

\item If $[\mathfrak{a}]=[\mathfrak{b}]$ then $\mathfrak{a} \ast Z_L \cong \mathfrak{b} \ast Z_L$.
If $\mathfrak{a}$ and $\mathfrak{b}$ represent the same elements in the class group, there exists
$\alpha \in F^\ast$ such that $\mathfrak{a}=\alpha \mathfrak{b}$. Hence $\mathfrak{a}L=\alpha \mathfrak{b}L$, and
$Z_{\mathfrak{a}L}$ and $Z_{\mathfrak{b}L}$ are isomorphic.

\item The action is simple.
If $Z_{\mathfrak{a}L}$ and $Z_{\mathfrak{b}L}$ are isomorphic, there exists $\alpha \in \B{R}$
such that 
\begin{equation}\label{toss}
\mathfrak{a}L=\alpha \mathfrak{b}L.
\end{equation}  
Recall that $L=\lambda \mathfrak{c}$ for some $\mathfrak{c} \in I_F$, and multiply both sides of \eqref{toss}
by the pseudolattice $\lambda^{-1} \mathfrak{b}^{-1}\mathfrak{c}^{-1}$.
This gives $\mathfrak{ab}^{-1} =\alpha \mathcal{O}_F$.
The left hand side is contained in $F$, so we have $\alpha \in F$. Hence $[\mathfrak{a}]=[\mathfrak{b}]$.

\item The action is transitive.  Let $Z_M$ be a quantum torus with endomorphism ring isomorphic to $\mathcal{O}_F$.
Then $M=\mu \mathfrak{d}$ for some $\mu \in \B{R}^\ast$, $\mathfrak{d} \in I_F$.  Put $\mathfrak{a}:=\mathfrak{dc}^{-1}$.
Then $\mathfrak{a} L=\lambda \mathfrak{ac}=\lambda \mathfrak{d}=(\lambda/\mu)M$. Hence
$\mathfrak{a} \ast Z_L \cong Z_M$.
\end{itemize}
\end{proof}

As a corollary we obtain Theorem \ref{C}:

\begin{proof}[Proof of Theorem \ref{C}]
The reciprocity map of class field theory gives us an isomorphism 
$$C(F) \cong \textrm{Gal}(H_F/F).$$
The result follows immediately from Theorem \ref{action}.
\end{proof}

We also obtain the following corollary:

\begin{cor}
$$\left|C(F) \right|=\left| \mathcal{QT}(\mathcal{O}_F)\right|.$$
\end{cor}

\noindent
{\bf Remarks: } 
\begin{enumerate}
\item Via \eqref{narrowsurjection} the narrow class group acts on $\mathcal{QT}(\mathcal{O}_F)$.
This is transitive, but is only faithful when $C(F)^+=C(F)$.  This occurs precisely when both infinite
primes of $F$ are unramified in the narrow ray class field of $F$. Class Field Theory gives us an
equation for the size of $C(F)^+$:
$$\left| C(F)^+ \right|=\frac{4 h_F }{[\mathcal{O}_F^\ast: F^+ \cap \mathcal{O}_F^\ast]},$$
where $h_F$ is the size of the class group of $F$.\\

Hence the action of $C(F)^+$ on $\mathcal{QT}(\mathcal{O}_F)$ is faithful precisely when
$$[\mathcal{O}_F^\ast: F^+ \cap \mathcal{O}_F^\ast]=4.$$
\item Let $E_{\Lambda}$ be an elliptic curve corresponding to a complex lattice $\Lambda$ by the Uniformization Theorem, and let $\sigma$ be an automorphism of $\B{C}$.  We have a natural action of $\sigma$ on $E_{\Lambda}$ by letting $\sigma$ act on the coefficients of the equation for $E_{\Lambda}$.  Analogous to the proof of Theorem \ref{action}, we have an action of the group of fractional ideals in $K$ on the set of elliptic curves $E$ with $\textrm{End}(E) \simeq \mathcal{O}_K$, where $K$ is some fixed quadratic imaginary field.  This action is defined and denoted by
$$(\mathfrak{a},E_\Lambda)\mapsto \mathfrak{a} \ast E_\Lambda:=E_{\mathfrak{a}\Lambda}.$$  If $\mathfrak{a}$ is a fractional ideal of $K$, then the reciprocity map supplies a homomorphism 
\begin{equation}\label{cg}
\psi_{H_K/K}: C(K) \rightarrow \textrm{Gal}(H_K/K)
\end{equation}
where $H_K$ is the Hilbert class field of $K$.  It is a fundamental result in the theory of Complex Multiplication that the following identity holds \cite{SilvermanII}:
$$ 
\mathfrak{a}^{-1} \ast E_{\Lambda} = E_{\Lambda}^{\psi_{H_K/K}([\mathfrak{a}])}. 
$$  
Since \eqref{cg} exhibits an isomorphism between the class group and the Galois group of the Hilbert class field of $K$ over $K$, this result is strongly linked to the following:
\begin{prop}
Let $E$ be an elliptic curve with Complex Multiplication by an order in an imaginary quadratic field $K$.  Then 
there exists an elliptic curve $E'$ such that $E$ and $E'$ are isomorphic, and $E'$ is defined over the Hilbert class field of $K$.
\end{prop}
By Theorem \ref{C} we are able to describe an action of $\textrm{Gal}(H_F/F)$ on isomorphism classes of Quantum Tori with Real Multiplication. There is no reason {\it a priori} why we should be able to do this.  The objects $Z_L$ are purely analytic constructions, associated to which there is no natural algebraic object for the automorphisms to act upon.  However, this simple result shows that Quantum Tori with Real Multiplication do possess algebraic characteristics.  Moreover, with respect to the algebraic property highlighted by this result, this suggests the existence of a virtual ``field of definition'' over which Quantum Tori are defined, and that in the case of Theorem {C} this field is equal to $H_F$.
\end{enumerate}

{\it Acknowledgements. }
The author would like to thank Boris Zilber and Nikolaos Diamantis for their encouragement and suggestions, and also Anand Pillay for his comments and advice.


\begin{thebibliography}{10}

\bibitem{Armstrong}
{\sc Armstrong, M.~A.}
\newblock {\em Basic Topology (Undergraduate Texts in Mathematics)}.
\newblock Springer, May 1997.

\bibitem{Connes}
{\sc Connes, A.}
\newblock {\em Noncommutative Geometry}.
\newblock {Academic Press}, November 1994.

\bibitem{BarnesII}
{\sc E.Barnes}.
\newblock The theory of the double gamma funcion.
\newblock {\em Philisophical Transactions of the Royal Society (A) 196\/}
  (1901), 265--388.

\bibitem{FesenkoII}
{\sc I.Fesenko}.
\newblock Model theory guidance in number theory?
\newblock To be published in Model Theory with Applications to Algebra and
  Analysis, LMS Lecture Note Series, CUP, 2007.

\bibitem{FesenkoI}
{\sc I.Fesenko}.
\newblock Several nonstandard remarks.
\newblock To be published in AMS/IP studies in Advances Mathmatics, AMS 2006.

\bibitem{Iwasawa}
{\sc Iwasawa, K.}
\newblock {\em Local Class Field Theory (Oxford Mathematical Monographs)}.
\newblock {Oxford University Press, USA}, September 1986.

\bibitem{Manin}
{\sc Manin, Y.~I.}
\newblock Real multiplication and noncommutative geometry (ein {A}lterstraum).
\newblock In {\em The legacy of Niels Henrik Abel}. Springer, Berlin, 2004,
  pp.~685--727.

\bibitem{Neukirch}
{\sc Neukirch, J.}
\newblock {\em Algebraic Number Theory (Grundlehren der mathematischen
  Wissenschaften)}.
\newblock Springer, June 1999.

\bibitem{ShimuraI}
{\sc Shimura, G.}
\newblock {\em Introduction to Arithmetic Theory of Automorphic Functions}.
\newblock {Princeton University Press}, August 1971.

\bibitem{ShimuraII}
{\sc Shimura, G.}
\newblock Class fields over real quadratic fields and {H}ecke operators.
\newblock {\em Ann. of Math. (2) 95\/} (1972), 130--190.

\bibitem{ShintaniIII}
{\sc Shintani, T.}
\newblock On certain ray class invariants of real quadratic fields.
\newblock {\em J. Math. Soc. Japan 30}, 1 (1978), 139--167.

\bibitem{SilvermanII}
{\sc Silverman, J.~H.}
\newblock {\em Advanced Topics in the Arithmetic of Elliptic Curves (Graduate
  Texts in Mathematics)}.
\newblock Springer, January 1994.

\end{thebibliography}
\end{document}